\documentclass{amsart}
\usepackage{amssymb, amsmath, amsfonts, amscd}
\usepackage{hyperref}

\input xypic

\theoremstyle{plain}
\newtheorem{theorem}{Theorem}[section]
\newtheorem{corollary}[theorem]{Corollary}
\newtheorem{lemma}[theorem]{Lemma}
\newtheorem{proposition}[theorem]{Proposition}

\newtheorem{example}[theorem]{Example}

\theoremstyle{definition}
\newtheorem{definition}[theorem]{Definition}

\theoremstyle{remark}

\numberwithin{equation}{theorem}

\renewcommand{\L}{\operatorname{L}}

\newcommand{\tp}{\tilde{\pi}}

\renewcommand{\O}{\mathcal{O} }

\newcommand{\Hom}{\operatorname{Hom}}

\newcommand{\Der}{\operatorname{Der} }

\newcommand{\End}{\operatorname{End} }
\newcommand{\Spec}{\operatorname{Spec} }
\renewcommand{\P}{\operatorname{P} }
\renewcommand{\H}{\operatorname{H} }
\newcommand{\D}{\operatorname{D^1(A,f)} }
\newcommand{\Dg}{\operatorname{D^1(A,g)}  }

\newcommand{\U}{\operatorname{U}} 
\newcommand{\K}{\operatorname{K} }

\newcommand{\Diff}{\operatorname{Diff}}

\newcommand{\Z}{\operatorname{Z} }
\newcommand{\tL}{\tilde{L} }
\newcommand{\ta}{\tilde{\alpha} }
\newcommand{\sym}{\operatorname{Sym} }

\newcommand{\LR}{\underline{\text{LR}(A/k)}}
\newcommand{\ATiA}{\underline{\operatorname{D^1(A,f)}-\text{Lie}} }
\newcommand{\Conn}{\underline{\operatorname{Conn}} }
\newcommand{\Mod}{\underline{\operatorname{Mod} }}
\newcommand{\Ext}{\operatorname{Ext}}
\newcommand{\CH}{\operatorname{CH}}
\newcommand{\Dl}{\operatorname{D} }

\newcommand{\C}{\operatorname{C} }

\newcommand{\Pic}{ \operatorname{Pic}  }
\newcommand{\Log}{ \operatorname{Log}  }
\newcommand{\Br}{ \operatorname{Br}  }

\newcommand{\fal}{f^{\alpha} }
\newcommand{\gal}{g^{\alpha} }

\newcommand{\ou}{\overline{u}}
\newcommand{\ov}{\overline{v}}
\newcommand{\Uo}{\operatorname{U}^{\otimes} }
\newcommand{\Ur}{\operatorname{U}^{\rho} }

\begin{document}

\title{Classification of Lie algebras of differential operators}

%\author{Helge \"{O}ystein Maakestad}

\author{Helge \"{O}ystein Maakestad \\ email \href{mailto:h_maakestad@hotmail.com}{\text{h\_maakestad@hotmail.com}} }

\email{\text{h\_maakestad@hotmail.com}}

\keywords{D-Lie algebra, canonical quotient, Lie-Rinehart algebra, connection, classification, non-abelian extension}

\thanks{}

\subjclass{}

\date{January 2019}

\begin{abstract} In a previous paper we introduced the notion of a $\Dl$-Lie algebra $\tL$. A $\Dl$-Lie algebra $\tL$ is an $A/k$-Lie-Rinehart algebra with a right $A$-module structure
and a canonical central element $D$ satisfying several conditions. We used this notion to define the \emph{universal enveloping algebra} of the category of $\tL$-connections
and to define the cohomology and homology of an arbitrary connection. In this note we introduce the canonical quotient $L$ of a $\Dl$-Lie algebra $\tL$ and use this to classify
$\Dl$-Lie algebras where $L$ is projective as $A$-module.  We define for any 2-cocycle $f\in \Z^2(\Der_k(A),A)$ a functor $F_{\fal}(-)$ from the category of 
$A/k$-Lie-Rinehart algebras to the category of $\Dl$-Lie algebras and classify $\Dl$-Lie algebras with projective canoncial quotient using the functor $F_{\fal}(-)$. We prove a similar classification for
non-abelian extensions of $\Dl$-Lie algebras. We moreover classify maps of $\Dl$-Lie algebras and $\tL$-connections $(E,\rho)$ in the case when the canonical quotient $L$ of $\tL$ is projective as $A$-module. 
\end{abstract}

\maketitle

\tableofcontents

\section{Introduction}

%If $X$ is a complex projective manifold and $E$ a finite rank holomorphic vector bundle on $X$ with a flat holomorphic connection $\nabla$, we may define the characteristic variety $\SS(E)$ of $E$.
%The variety $\SS(E)$ is a sub variety of $\mathbb{V}(\Omega^1_X)$ - the cotangent bundle of $X$, and $\SS(E)$ is a coisotropic variety with respect to the natural Poisson product on $\mathbb{V}(\Omega^1_X)$.
%The aim of this paper is to generalize this construction to the algebraic situation and the case when the connection $\nabla$ is non-flat. We do this in the language of $\Dl$-Lie algebras 
%and connections on $\Dl$-Lie algebras. 

A $\Dl$-Lie algebra $\tL$ is a refinement of the notion of an $A/k$-Lie-Rinehart algebra: It is a $k$-Lie algebra and left $\P^1_{A/k}$-module, where $\P^1_{A/k}$ is the first
order module of principal parts of $A/k$. The Lie-product on $\tL$ satisfies the anchor condition 
\[  [u,cv]=c[u,v]+\tp(u)(c)v \]
for all $u,v\in \tL$ and $c\in A$. Here $\tp: \tL \rightarrow \Der_k(A)$ is a map of $k$-Lie algebras and left $\P$-modules. Hence if we view $\tL$ as a left $A$-module and $k$-Lie algebra, it follows
the pair $(\tL, \tp)$ is an $A/k$-Lie-Rinehart algebra. There is a canonical central element  $D\in \tL$ with the property that 
\[ uc=cu+\tp(u)(c)D \]
holds for $u\in \tL$ and $c\in A$. We may define the notion of a connection $\rho: \tL \rightarrow \Diff^1(E)$ where $E$ is a left $A$-module. In the papers \cite{maa1} and \cite{maa2} general properties
of the category of $\Dl$-Lie algebras are introduced and studied: Universal enveloping algebras of $\Dl$-Lie algebras, cohomology and homology of connections on $\Dl$-Lie algebras. 
Let $\underline{\operatorname{D-Lie}}$ denote the category of $\Dl$-Lie algebras and morphisms. In \cite{maa1} Theorem 3.5 and 3.18 we prove the following result:

\begin{theorem}\label{main1}  There are covariant functors 
\begin{align}
&\Uo: \underline{\operatorname{D-Lie}} \rightarrow \underline{Rings} \\
&\Ur: \underline{\operatorname{D-Lie}} \rightarrow \underline{Rings}.
\end{align}
with the following property: For any $\Dl$-Lie algebra $\tL$ there are exact equivalences of categories
\begin{align}
& F_1:\Mod(\tL, Id) \cong \Mod(\Uo(\tL)) \\
& F_2: \Conn(\tL, Id) \cong \Mod(\Ur(\tL))
\end{align}
with the property that $F_1$ and $F_2$ preserves injective and projective objects.
\end{theorem}
We used the associative unital rings $\Uo(\tL)$ and $\Ur(\tL)$ to define the cohomology and homology of an arbitrary $\tL$-connection $(E,\rho)$.  By Lemma 3.19 in \cite{maa1} there is for any 2-cocycle $f \in \H^2(\Der_k(A),A)$ an exact equivalence of categories
\[ \psi_f: \Conn(L) \cong \Mod(\Uo(L(\fal)) \]
preserving injective and projective objetcs. Here $L$ is an $A/k$-Lie-Rinehart algebra and $\Conn(L)$ is the category of $L$-connections and morphisms of connections. Hence we may use the associative unital ring $\Uo(L(\fal))$ to define the cohomology and homology of any $L$-connection, flat or non-flat. Using $\psi_f$ we get for any pair of $L$-connections $V,W$ and any integer $i \geq 0$ isomorphisms
\[   \Ext^i_{\Uo(L(\fal))}(\psi_f(V), \psi_f(W)) \cong   \Ext^i_{\Conn(L)}(V,W) .\]
Hence the associative ring $\Uo(L(\fal))$ may be used to calculate the "true" $\Ext$-groups of $V$ and $W$. Theorem \ref{main1} was one of the main reasons for the introduction of the notion $\Dl$-Lie algebra.

Given an almost commutative ring $U$ with filtration $U_i$ where $U$ is generated by $U_1$, there is a pre-$\Dl$-Lie algebra $\tL$ and a 2-sided ideal $J\subseteq \Uo(\tL)$ and an isomorphism
\[ \Uo_J(\tL) \cong U \]
of filtered associative rings. Hence there is an exact equivalence of categories
\[ \Mod(U) \cong \Mod(\tL, Id, J) \]
preserving injective and projective objects, where $\Mod(U)$ is the category of left $U$-modules and $\Mod(\tL, Id, J)$is the category of $\tL$-connections with $J$-curvature equal to zero. 
Most rings of differential operators appearing in the litterature
are almost commutative and generated by $U_1$. As special cases it follows the ring $\Diff(\mathcal{L})$ of differential operators on a linebundle $\mathcal{L}$ on $A$ 
and  the generalized universal enveloping algebra  $\U(A,L,f)$ where $L$ is an $A/k$-Lie-Rinehart algebra and $f$ is a 2-cocycle  on $L$ with values in $A$,
 may be realized as $\Uo_J(\tL)$ (see \cite{maa1}, Example 3.35). The rings $\Uo(\tL)$ and $\Uo_J(\tL)$ have canonical filtrations $\Uo(\tL)^i, \Uo_J(\tL)^i$ and the associated graded rings are commutative.
Let $\tL^*:=\Uo(\tL)^1/\Uo(\tL)^0$. There is a canonical surjective map of graded $A$-algebras
\[ \rho: \sym_A^*(\tL^*) \rightarrow Gr(\Uo(\tL)) .\]
The map $\rho$ is not always an isomorphism in the case when $\tL$ is a projective $A$-module. Hence the PBW-theorem does not always hold for $\Uo(\tL)$. If $(U,U_i)$ is an almost commutative ring with $L:=U_1/U_0$ a projective $A$-module and the canonical map $\eta:\sym_A^*(L) \rightarrow Gr(U)$ an isomorphism, it follows there is an isomorphism of associative rings $U \cong \U(A,L,f)$ where $L:=U_1/U_0$ and
$f$ a 2-cocycle for $L$ with values in $A$. Hence a left $U$-module $A$ corresponds to a connection $\nabla: L\rightarrow \End_k(E)$ of curvature type $f$. It follows the k'th Chern class of $E$ is determined
by the first Chern class. A similar PBW-theorem for the universal ring $\Uo(\tL)$ would put similar conditions on \emph{all Chern classes} of all finite rank projective $A$-modules. Hence the study
of the associative ring $\Uo(\tL)$ and the PBW-theorem for $\Uo(\tL)$ will have applications in the study of the Chern classes and Chern character of $A$.

Note: The category $\Conn(L)$ of connections on an $A/k$-Lie-Rinehart algebra $L$ is a small abelian category, hence the Freyd-Mitchell full embedding theorem gives an equivalence $\phi$  of $\Conn(L)$ with
a sub category of $\Mod(R)$ where $R$ is an associative ring. The equivalence $\phi$ does not preserve injective and projective objects, hence we cannot use $\phi$ to define the cohomology and homology
of a connection. Theorem \ref{main1} gives a geometric construction of cohomology and homology groups of any $\tL$- and $L$-connection where $\tL$ is any $\Dl$-Lie algebra and
$L$ is any $A/k$-Lie-Rinehart algebra. The construction is functorial and can be done for any scheme and any  sheaf of Lie-Rinehart algebras.

The aim of this paper is to classify maps of $\Dl$-Lie algebras. We also classify $\Dl$-Lie algebras and non-abelian extensions of $\Dl$-Lie algebras with projective canonical quotient in terms of a functor $F_{\fal}$ introduced in \cite{maa1}.

Let $(L,\alpha)$ and $(L',\alpha')$ be $A/k$-Lie-Rinehart algebras and let $f,g\in \Z^2(\Der_k(A),A)$ be two 2-cocycles. Let $F_g(L):=L(\gal)$ and $F_f(L'):=L'(f^{\alpha'})$ be the 
$\Dl$-Lie algebras introduced in \cite{maa1}, Theorem 2.7. In Theorem \ref{mapclass} we prove the following:

\begin{theorem} There is a map of $\Dl$-Lie algebras $\phi: L(\gal)\rightarrow L'(f^{\alpha'} )$ if and only if $\overline{\gal}=\overline{\fal}$ in $\alpha^*\H^2(\Der_k(A),A)$. In this case 
there is an equality between the set of maps of $\Dl$-Lie algebras $\phi: L(\gal)\rightarrow L'(f^{\alpha'} )$ and the set of maps of $A/k$-Lie-Rinehart algebras $\psi_2:L \rightarrow L'$.
\end{theorem}

In Theorem \ref{Dclassification} we prove the following:

\begin{theorem} Let $(\tL, \ta,\tp, [,],D)$ be a $\Dl$-Lie algebra and let $(L,\alpha)$ be the canonical quotient $A/k$-Lie-Rinehart algebra of $\tL$. Assume $L$ is projective as left $A$-module. 
There is an isomorphism $\tL \cong L(\fal)$ as $\Dl$-Lie algebras where $L(\fal):=F_{\fal}(L)$ and $F_{\fal}$ is the functor from Example \ref{nontrivial}. Hence $\tL$ is uniquely determined up to isomorphism
by the canonical quotient $(L,\alpha)$ and the 2-cocycle $f\in \H^2(\Der_k(A),A)$.
\end{theorem}

In Corollary \ref{mapclasscorr} we classify maps of $\Dl$-Lie algebras between two $\Dl$-Lie algebras $\tL_1$ and $\tL_2$.

In Corollary \ref{nonabclass} we get a classification of non-abelian extensions of $\Dl$-Lie algebras in terms of the functor $F_{\fal}$.

% and the aim of this paper is to define and study the characteristic variety 
%$\SS(E,\nabla)$ of a connection $\nabla$ that is non-flat using non-abelian extensions of $\Dl$-Lie algebras.

In the paper we introduce the category of $\Dl$-Lie algebras, connections on $\Dl$-Lie algebras using the module of principal parts. This gives an equivalent definition of a $\Dl$-Lie algebra to the one introduced in \cite{maa1}. We prove in Theorem \ref{mthm} that for any $\tL$-connection $(E,\rho)$ there is an extension of $\Dl$-Lie algebras
\begin{align}
&\label{nonabel} 0 \rightarrow \End_A(E) \rightarrow \End(\tL, E) \rightarrow \tL \rightarrow 0 
\end{align}
with the  property that \ref{nonabel} splits in the category of $\Dl$-Lie algebras if and only if $E$ has a flat $\tL$-connection $\rho':\tL \rightarrow \Diff^1(E)$. Theorem \ref{mthm} generalize a similar result proved
in \cite{maa2} for non-abelian extensions of $A/k$-Lie-Rinehart algebras.

% We introduce the canonical quotient $(L,\alpha)$ of a $\Dl$-Lie algebra $\tL$ and prove in Theorem \ref{Dclassification} there is an exact sequence of $A/k$-Lie-Rinehart algebras
%\[ 0 \rightarrow J \rightarrow \tL \rightarrow L \rightarrow 0 \]
%with the property that if $L$ is projective as left $A$-module there is a 2-cocycle $f\in \Z^2(\Der_k(A),A)$ and an isomorphism $\tL \cong L(\fal)$ where $L(\fal):=F_{\fal}(L)$, and where 
%$F_{\fal}$ is the functor defined and studied in \cite{maa1}. Hence when $L$ is projective it follows $\tL$ is uniquely determined by $L$ and the cohomology class $\overline{f}\in \H^2(\Der_k(A),A)$.
%Hence Theorem \ref{Dclassification} classifies $\Dl$-Lie algebras in terms of 2-cocycles and $A/k$-Lie-Rinehart algebras. In Corollary \ref{nonabclass} we get a similar classification of non-abelian extensions of $\Dl$-Lie %algebras in terms of the functor $F_{\fal}$.

Using the $\P$-module structure on the $\Dl$-Lie algebra $\tL$ we define in Definition \ref{correspond} the correspondence
$Z(\rho, \sigma_1,\sigma_2)$ associated to a connection and degeneracy locies $\sigma_1$ and $\sigma_2$. The correspondence $Z(\rho, \sigma_1,\sigma_2)$ induce an endomorphism
\[ I(\rho,\sigma_1,\sigma_2): \CH^*(X)\rightarrow \CH^*(X) \]
of the Chow-group $\CH^*(X)$. The correspondence $Z(\rho, \sigma_1,\sigma_2)$ and Chow-operator $I(\rho,\sigma_1,\sigma_2)$ depend in a non-trivial way on the $\P$-module structure of $\tL$ and $\P$-linearity of the $\tL$-connection $\rho$ and cannot be defined
for an ordinary $L$-connection where $(L,\alpha)$ is an $A/k$-Lie-Rinehart algebra (see Example \ref{correspondence}). There is no non-trivial right $A$-module structure on $L$.

We also prove in Theorem \ref{LconnLfalconn} the following result relating the category of connections on an $A/k$-Lie-Rinehart algebra $L$ to the category of connections on the $\Dl$-Lie algebra $L(\fal)$:

\begin{theorem}\label{LconnLfalconn}  Let $(L,\alpha)$ be an $A/k$-Lie-Rinehart algebra and let $f\in \Z^2(\Der_k(A),A)$ be a 2-cocycle. There is an equivalence of categories
\[ C_f: \Conn(L,\End) \rightarrow \Conn(L(\fal)) \]
from the category $\Conn(L,\End)$ of $(L,\psi)$-connections $\nabla$, to the category $\Conn(L(\fal))$ of $L(\fal)$-connections $\rho$.
Let $u:=az+x, v:=bz+y\in L(\fal)$ and let $e\in E$. Let $(E,\nabla)$ be an $L$-connection and let $\rho_{\nabla}:=C_f(\nabla)$. The following holds:
\[ R_{\rho_{\nabla}}(u,v)(e) =R_{\nabla \circ i}(x,y)(e)-\fal(u,v)e.\]
For any $L(\fal)$-connection $(E,\rho)$ there is an $(L,\psi)$-connection $(E,\nabla)$ with $C_f(E,\nabla)=(E,\rho)$.
\end{theorem}

Hence there is for any 2-cocycle $f\in \Z^2(\Der_k(A),A)$ a functorial way to define an $L(\fal)$-connection $C_f(E,\nabla)$ from an $L$-connection $(L,\alpha)$.
As a Corollary we are able to classify $\tL$-connections on $\Dl$-Lie algebras with projective canonical quotient. In Corollary \ref{DclassCorr} we prove the following:

\begin{corollary} Let $(\tL,\ta,\tp,[,],D)$ be a $\Dl$-Lie algebra and let $(E,\rho)$ be an $\tL$-connection. Assume the canonical quotient $L$ of $\tL$
is a projective $A$-module. It follows any $\tL$-connection $(E,\rho)$ is on the form $C_f(E,\nabla)$ where $f\in \Z^2(\Der_k(A),A)$ and $(E,\nabla)$ is an $(L, \psi)$-connection for $\psi\in \End_A(E)$.
If $\rho(D)=I$ it follows we may choose $(E,\nabla)$ to be an $L$-connection.
\end{corollary}

\section{Classification of morphisms of $\Dl$-Lie algebras}

Let in the following $A$ be a fixed commutative $k$-algebra where $k$ is a fixed commutative unital ring. Let $f,g\in \Z^2(\Der_k(A),A)$ bw two 2-cocycles and consider the two
$k$-Lie algebras $\D$ and $\Dg$. Let $\P:=A\otimes_k A/I^2$ be the first order module of principal parts of $A/k$. The $k$-Lie algebra $\D$ has the following left and right $A$-module structure:
\[ c(a,x):=(ca,cx) \]
and
\[ (a,x)c:=(ac+x(c),cx) \]
for $(a,x)\in \D$ and $c\in A$. Let $z_f:=(1,0)\in \Dl^1(A,f)$. It follows $u:=(a,x)$ that  $uc=cu+x(c)z_f=cu+\pi_f(u)(c)z_f$. Similarly for $\Dg$. Let $d:A \rightarrow A\otimes_k A$ be defined by
$dc:=1\otimes c-c\otimes 1$. It follows $dc(u)=uc-cu=\pi_f(u)(c)z$. One checks that $zc=cz$ hence $dbdc(u)=db(\pi_f(u)(c)z=0$ hence $\D$ is annihilated by $I^2\subseteq A\otimes_k A$
and it follows $\D$ is a left $\P$-module. 
The two projection maps $\pi_f:\D \rightarrow \Der_k(A)$ and $\pi_g: \Dg \rightarrow \Der_k(A)$ and 
$\pi_f,\pi_g$ are maps of $\P$-modules and $k$-Lie algeras. 
We want to classify maps
\[ \phi: \Dg  \rightarrow \D \]
of $\P$-modules and $k$-Lie algebras such that 

\begin{align}
&\label{mapcond} \pi_f \circ \phi = \pi_g \text{ and  }\phi(z_g)=z_f.
\end{align}

\begin{lemma} Let $\phi: \Dg  \rightarrow \D$ be a map of $\P$-modules and $k$-Lie algebras satisfying condition \ref{mapcond}. It follows $\phi(a,x)=(a+\phi_1(x),x)$ where
 $\phi_1\in \C^1(\Der_k(A),A)$, $ g=f+d^1(\phi_1)$ and $d^1: \C^1(\Der_k(A),A)\rightarrow \C^2(\Der_k(A),A)$ is the first differential in the Lie-Rinehart complex of $\Der_k(A)$. Hence the map $\phi$ exists if and only if
$\overline{g}=\overline{f}\in \H^2(\Der_k(A),A)$. If the map $\phi$ exists it is an isomorphism of $\P$-modules and $k$-Lie algebras with inverse $\psi:\D  \rightarrow \Dg$ defined by
$\psi(a,x):=(a-\phi_1(x),x)$.
\end{lemma}
\begin{proof} Let $\phi : \Dg \rightarrow \D$ be a map of $\P$-modules and $k$-Lie algebras satisfying \ref{mapcond} It follows for any $(a,x)\in \Dg$ we get
\begin{align}
&\label{map} \phi(a,x)=\phi(a,0)+\phi(0,x)=az_f+(\phi_1(x), x)=(a+\phi_1(x),x) 
\end{align}
with $\phi_1 \in \C^1(\Der_k(A),A)$. One checks the map $\phi$ is a map of $k$-Lie algebras if and only if $g=f+d^1(\phi_1)$. The map $\phi$ in \ref{map} is always a map of $\P$-modules.
Hence the map $\phi$ exists if and only if $\overline{g}=\overline{f}\in \H^2(\Der_k(A),A)$. If $\phi$ exists, an inverse $\psi$ is given by
$\psi(a,x):=(a-\phi_1(x),x)$. The Lemma follows.
\end{proof}

Assume $(L,\alpha)$ is an $A/k$-Lie-Rinehart algebra and let $f,g\in \Z^2(\Der_k(A),A)$ be two 2-cocycles and let $\fal,  \gal\in \Z^2(L,A)$ be the pull-back cocycles. Let $F_g(L):=L(\gal):=Az_g\oplus L$ 
be the $A/k$-Lie-Rinehart algebra defined in \cite{maa1}. It has the following $k$-Lie product:
\[ [(a,x),(b,y)]:=(\alpha(x)(b)-\alpha(y)(a)+\gal(x,y),[x,y])\]
for $(a,x),(b,y)\in L(\gal)$. There is a map 
\[ \pi_g: L(\gal) \rightarrow \Der_k(A) \]
defined by
\[ \pi_g(a,x):=\alpha(x).\]
The pair $(L(\gal), \pi_g)$ is an $A/k$-Lie-Rinehart algebra. 

Note: The map $\alpha$ induce a map of cohomology groups
\[ \alpha^*: \H^2(\Der_k(A),A) \rightarrow \H^2(L,A).\]
Let $Im(\alpha^*):=\alpha^*\H^2(\Der_k(A),A)$ denote the image of the mape $\alpha^*$ in $\H^2(L,A)$.

\begin{lemma} \label{map1}There is a map $\alpha_g: L(\gal) \rightarrow \D$ of $\P$-modules and $k$-Lie algebras with $\pi \circ \alpha_g =\pi_g$ if and only if there is an element $\alpha_1\in \C^1(L,A)$
with $\gal =\fal+d^1(\alpha_1)$. The map $\alpha_g$ is defined as follows: $\alpha_g(a,x):=(a+\alpha_1(x),\alpha(x))\in \D$. Hence there exists a 5-tuple $(L(\gal), \alpha_g, \pi_g,[,],z_g)$ which is a
$\Dl$-Lie algebra if and only if there is an equality of cohomology classes 
\[ \overline{\gal}=\overline{\fal}\in \alpha^*\H^2(\Der_k(A),A) \subseteq \H^2(L,A).\]
\end{lemma}
\begin{proof} By definition the map $\alpha_g$ must look as follows:
\[ \alpha_g(a,x)=(a+\alpha_1(x), \alpha(x) )\]
where $(a,x)\in L(\gal)$ and $\alpha_1\in \C^1(L,A)$. One checks the map $\alpha_g$ is a map of left $\P$-modules. The map $\alpha_g$ is a map of $k$-Lie algebras
if and only if $\gal=\fal +d^1(\alpha_1)$. Hence the 5-tuple $(L(\gal), \alpha_g, \pi_g,[,],z_g)$ is a $\Dl$-Lie algebra over $\D$ if and only if $\overline{\gal}=\overline{\fal}$ in 
$\alpha^*\H^2(\Der_k(A),A)$. The Lemma follows.
\end{proof}

We may now classify maps between arbitrary $\Dl$-Lie algebras. Let $(L,\alpha)$ and $(L',\alpha')$ be $A/k$-Lie-Rinehart algebras and let $f,g\in \Z^2(\Der_k(A),A)$ with $\gal =\fal+d^1(\alpha_1)$
for $\alpha_1\in \C^1(L,A)$. Hence there is by Lemma \ref{map1} a structure as $\Dl$-Lie algebra $\alpha_g: L(\gal) \rightarrow \D$ given by the map $\alpha_g(a,x):=(a+\alpha_1(x),\alpha(x))$.
The 5-tuple $(L(\gal), \alpha_g,\pi_g,[,],z_g)$ is a $\Dl$-Lie algebra over $\D$ with $\pi_g(a,x):=\alpha(x)\in \Der_k(A)$.

Let $(L,\alpha)$ and $(L',\alpha')$ be $A/k$-Lie-Rinehart algebras and let $f,g\in \Z^2(\Der_k(A),A)$ be two 2-cocycles. Let $F_g(L):=L(\gal)$ and $F_f(L'):=L'(f^{\alpha'})$ be the 
$\Dl$-Lie algebras introduced in \cite{maa1}, Theorem 2.7.

\begin{theorem} \label{mapclass}  There is a map of $\Dl$-Lie algebras $\phi: L(\gal)\rightarrow L'(f^{\alpha'} )$ if and only if $\overline{\gal}=\overline{\fal}$ in $\alpha^*\H^2(\Der_k(A),A)$. In this case 
there is an equality between the set of maps of $\Dl$-Lie algebras $\phi: L(\gal)\rightarrow L'(f^{\alpha'} )$ and the set of maps of $A/k$-Lie-Rinehart algebras $\psi_2:L \rightarrow L'$.
\end{theorem}
\begin{proof} If there is an equality $\overline{\gal}=\overline{\fal}$ it follows there is a map
\[ \alpha_g: L(\gal) \rightarrow \D \]
defined by
\[ \alpha_g(a,x)=(a+\alpha_1(x), \alpha(x)).\]
The map $\phi: L(\gal) \rightarrow L'(f^{\alpha'})$ must look as follows:
\[ \phi(a,x)=(a+\phi_1(x),\phi_2(x)) \]
where $\phi_1: L \rightarrow A$ and $\phi_2: L \rightarrow L'$. Since $\alpha'_f \circ \phi = \alpha_g$ it follows $\phi_1=\alpha_1$ and $\alpha' \circ \phi_2 = \alpha$. Hence $\phi_2$ is a map
of $A/k$-Lie-Rinehart algebras. One checks the map $\phi$ is $\P$-linear and a map of $k$-Lie algebras. Conversely, for an arbitrary map $\phi_2: L\rightarrow L'$ of 
$A/k$-Lie-Rinehart algebras it follows the map $\phi(a,x):=(a+\alpha_1(x), \phi_2(x))$ is a map of $\Dl$-Lie algebras. The Theorem follows.
\end{proof}

\section{Classification of $\Dl$-Lie algebras with projective canonical quotient}

In this section we define the notion of a $\Dl$-Lie algebra, the category of $\Dl$-Lie algebras and connections on $\Dl$-Lie algebras. 

The module $\Diff^1(A)$ of first order differential operators
on a commutative $k$-algebra $A$ is in a canonical way a $k$-Lie algebra and a left $\P$-module where $\P:=A\otimes_k A/I^2$ is the module of first order principal parts on $A/k$. 
By definition $\Diff^1(A):=\Hom_A(\P,A)$ and it follows $\Diff^1(A)$ is in a canonical way a left and right $A$-module. There is a canonical projection map $\pi: \Diff^1(A) \rightarrow \Der_k(A)$ which is a map of $k$-Lie algebras and left $\P$-modules. There is a canonical inclusion of left and right $A$-modules and $k$-Lie algebras $\Diff^1(A) \subseteq \End_k(A)$ and 
the Lie product $[,]$ on $\Diff^1(A)$ satisfies the following formula:
\begin{align}
&\label{a1}  [u,cv]=c[u,v]+\pi(u)(c)v 
\end{align}
for all $u,v\in \Diff^1(A)$ and $c\in A$. There is a canonical element $D\in \Diff^1(A)$ with the property that 
\begin{align}
&\label{a2} uc=cu+\pi(u)(c) D.
\end{align}
The element $D$ is central in $\Diff^1(A)$ and $\pi(D)=0$. A $\Dl$-Lie algebra $\tL$ is a generalization of $\Diff^1(A)$ and we may speak the category of $\Dl$-Lie algebras, extensions and non-abelian extensions 
of $\Dl$-Lie algebras, cohomology of $\Dl$-Lie algebras, connections etc. Equation \ref{a1} is the equation defining an $A/k$-Lie-Rinehart algebra and any $\Dl$-Lie algebra $\tL$ has an underlying $A/k$-Lie-Rinehart algebra. Hence we may view a $\Dl$-Lie algebra as a refinement of the notion of a $A/k$-Lie-Rinehart
algebra. A $\Dl$-Lie algebra is an $A/k$-Lie-Rinehart algebra with extra structure defined by the $\P$-module structure, the canonical element $D$ and Equations \ref{a1} and \ref{a2}.

Note: We may consider the $A/k$-Lie-Rinehart algebra $\Der_k(A)\subseteq \Diff^1(A)$ and $\Der_k(A)$ has a canonical structure as $k$-Lie algebra and left $A$-module. It has no-non-trivial right
$A$-module structure. To get a non-trivial right $A$-module structure, we must consider the abelian extension $\Diff^1(A)$ and this is one of the motivations for the  introduction of the notion
of a $\Dl$-Lie algebra. The $\P$-module structure on $\Diff^1(A)$ is canonical and is related to the notion of curvature of a connection. 

Using the $\P$-module structure on $\tL$ we define the correspondence
$Z(\rho, \sigma_1,\sigma_2)$ associated to a connection and degeneracy locies $\sigma_1$ and $\sigma_2$. The correspondence $Z(\rho, \sigma_1,\sigma_2)$ induce an endomorphism
\[ I(\rho,\sigma_1,\sigma_2): \CH^*(X)\rightarrow \CH^*(X) \]
of the Chow-group $\CH^*(X)$.

The main Theorem in this section is Theorem \ref{mthm} where we construct the non-abelian extension $\End(\tL,E)$ of any $\tL$-connection $(E,\rho)$ in $\Conn(\tL, Id)$,
 and prove the canonical projection map of $\Dl$-Lie algebras 
\[ p_E: \End(\tL,E) \rightarrow \tL \]
splits in the category of $\Dl$-Lie algebras if and only if $E$ has a flat $\tL$-connection in $\Conn(\tL,Id)$. In \cite{maa2} a similar result was proved for $A/k$-Lie-Rinehart algebras.

We also define the canonical quotient $A/k$-Lie-Rinehart algebra $(L,\pi_L)$ of a $\Dl$-Lie algebra $\tL$ and classify the set of $\Dl$-Lie algebras with canonical quotient $L$. It follows from Theorem
\ref{Dclassification} the $\Dl$-Lie algebra $\tL$ is uniquely determined by $(L, \pi_L)$ up to isomorphism when $L$ is projective as left $A$-module. Hence for any cohomology class
$c=\overline{f} \in \H^2(\Der_k(A),A)$ and any $A/k$-Lie-Rinehart algebra $(L,\alpha)$ with $L$ a projective $A$-module, there is a unique $\Dl$-Lie algebra $\ta: \tL \rightarrow \D$ with canonical quotient $(L,\pi_L)$, defined by $\tL:=F_{\fal}(L)$, where $F_{\fal}$ is the functor defined and studied in \cite{maa1}.

Let in the following $A$ be a fixed commutative unital $k$-algebra where $k$ is a commutative unital ring. Let $\Der_k(A)$ be the $k$-Lie algebra of $k$-linear derivations of $A$.
Let $f\in \Z^2(\Der_k(A),A)$ be a 2-cocycle and let $\D:=A\oplus \Der_k(A)$ with the following $k$-Lie algebra structure and $A\otimes_k A$-module structure:
Let $u:=(a,x), v:=(b,y)\in \D$ and define
\[ [u,v]:=(x(a)-y(b)+f(x,y),[x,y]) \in \D.\]
Define for any element $c\in A$
\[ cu:=(ca,cx) \]
and
\[ uc:=(ac+x(c),cx).\]
It follows $\D$ is a left $A\otimes_k A$-module. Define the map $\pi: \D\rightarrow \Der_k(A)$ by $\pi(u):=\pi(a,x)=x$. Define for $x\in \Der_k(A)$ $xc:=cx$. It follows $\Der_k(A)$ is a left  $A\otimes_k A$-module.
Let $\P:=A\otimes_k A/I^2$ where $I$ is the kernel of the multiplication map. It follows $\P$ is the first order module of principal parts of $A/k$. Let $d:A\rightarrow \P$ be the universal derivation.
Let $p(a):=1\otimes a$ and $q(a):= a\otimes 1$. It follows $d=p-q$. Let $D:= (1,0)\in \D$.

\begin{lemma} \label{origin}The $k$-Lie algebra $\D$ is in a canonical way a left $\P$-module. The map $\pi: \D\rightarrow \Der_k(A)$ is a map of $k$-Lie algebras and left $\P$-modules.
The following holds for all $u,v\in \D$ and $c\in A$:
\begin{align}
&\label{eq1}  [u,cv]=c[u,v]+\pi(u)(c)v  \\
&\label{eq2} dc.u=\pi(u)(c)D 
\end{align}
The element $D$ is a central element with $\pi(D)=0$.
\end{lemma}
\begin{proof} Since $\D$ is a left and right $A$-module with $(au)b=a(ub)$ for all $u\in \D$ and $a,b\in A$ it follows $\D$ is a left $A\otimes_k A$-module. One checks that for any element 
$w\in I^2$ it follows $wu=0$ hence $\D$ is a left $\P$-module. The map $\pi$ is left and right $A$-linear. It follows $\pi$  is a map of $\P$-modules and $k$-Lie algebras. One checks 
Equation \ref{eq1} and \ref{eq2} holds and the Lemma is proved.
\end{proof}

The notion of a $\Dl$-Lie algebra is a generalization of the $k$-Lie algebra and left $\P$-module $\D$ from Lemma \ref{origin}.

\begin{definition} \label{mainD} Let $f\in \Z^2(\Der_k(A),A)$ be a 2-cocycle and let $\D$ be the $k$-Lie algebra and $\P$-module defined in Lemma \ref{origin}.
A 5-tuple $(\tL, \ta, \tp, [,],D)$ is a \emph{$\Dl$-Lie algebra} if the following holds. $\tL$ is a $k$-Lie algebra and left $\P$-module. The element $\ta$ is a map
\[ \ta: \tL \rightarrow \D \]
of $\P$-modules and $k$-Lie algebras. The element $\tp$ is a map
\[ \tp: \tL \rightarrow \Der_k(A) \]
of left $\P$-modules and $k$-Lie algebras with $ \tp= \pi \circ \ta$. The $k$-Lie product satisfies 
\[ [u,av]=a[u,v]+\tp(u)(a)v \]
for all $u,v\in \tL$ and $a\in A$. The following holds for all $a\in A$:
\[ da.u=\tp(u)(a)D .\]
The element $D$ is a central element in $\tL$ with $\tp(D)=0$. Given two $\Dl$-Lie algebras $(\tL, \ta,\tp,[,],D)$ and $(\tL', \ta',\tp',[,],D')$.  A \emph{map of $\Dl$-Lie algebras} is a map
\[ \phi: \tL\rightarrow \tL' \]
of left $\P$-modules  and $k$-Lie algebras with $\phi(D)=D'$ and $\tp' = \tp \circ \phi$. Let \underline{$\Dl$-Lie} denote the category of $\Dl$-Lie algebras and maps. An \emph{ideal} $\tilde{I}$ 
in a $\Dl$-Lie algebra $\tL$ is a sub-$k$-Lie algebra and sub-$\P$-module $\tilde{I} \subseteq \tL$ such that $[\tL, \tilde{I}]\subseteq \tilde{I}$. A \emph{$\Dl$-ideal} $\tilde{I}$ is a 
sub-$k$-Lie algebra and sub-$\P$-module $\tilde{I} \subseteq \tL$ with $[\tL, \tilde{I} ] \subseteq \tilde{I}$ and with $D\notin \tilde{I}$. 
%A $\Dl$-Lie algebra $\tL$ is a \emph{trivial $\Dl$-Lie algebra} if $D=0$.
\end{definition}

Note: One checks Definition \ref{mainD} gives the same notion as Definition 2.3 in \cite{maa1}.

\begin{example} \label{nontrivial}  Non-trivial examples: $\Dl$-Lie algebras and $A/k$-Lie-Rinehart algebras.\end{example}

 It follows from Lemma \ref{origin} that the 5-tuple $(\D, id, \pi, [,],D)$ is a $\Dl$-Lie algebra for any 2-cocycle $f\in \Z^2(\Der_k(A),A)$. 

Given a $\Dl$-Lie algebra 
$(\tL, \ta, \tp,[,],D)$ it follows the left $A$-module $\tL$ and the map $\tp: \tL \rightarrow \Der_k(A)$ is an $A/k$-Lie-Rinehart algebra. 

Given any 2-cocycle $f\in \Z^2(\Der_k(A),A)$, let $(L,\alpha)$ be an $A/k$-Lie-Rinehart algebra and let $\fal \in \Z^2(L,A)$ be the pull back 2-cocycle. There is by Theorem 2.7 in \cite{maa1} a functor 
\[ F_{f^{\alpha} }: \LR  \rightarrow \underline{\Dl-\text{Lie}} \]
from the category \LR \text{ }of $A/k$-Lie-Rinehart algebras to the category \underline{$\Dl$-Lie} of $\Dl$-Lie algebras, hence it is easy to give non-trivial examples of $\Dl$-Lie algebras. When $\H^2(\Der_k(A),A)\neq 0$ we get for any $f\in \Z^2(\Der_k(A),A)$ a non-trivial functor $F_{\fal}$. Since $F_{\fal}(L)$ is independent of choice of representative for the class $c:=\overline{\fal}\in \H^2(L,A)$ the following holds: 
There is for any two representatives $f,f'$ for the class $c$ a canonical isomorphism $F_{f}(L) \cong F_{f'}(L)$ of extensions of $A/k$-Lie-Rinehart algebras.
Hence the two functors $F_f$ and $F_{f'}$ are equal up to isomorphism. Hence we get from Theorem 2.7 in \cite{maa1} a well defined functor $F_c$ for any cohomology class $c:=\overline{\fal}\in \H^2(L,A)$:
\[ F_c: \LR \rightarrow \underline{\Dl-\text{Lie}} .\]
When the class $c$ is the zero class it follows the underlying $A/k$-Lie-Rinehart algebra of $F_c(L):=Az\oplus L$ is the trivial abelian extension of $L$ by $A$. Hence the cohomology group $\H^2(L,A)$ parametrize
a large class of non-trivial functors $F_c$.

%: \LR \rightarrow  \underline{\Dl-\text{Lie}} $.

\begin{lemma} \label{lemma1} Let $\phi: \tL_1\rightarrow \tL_2$ be a map of $\Dl$-Lie algebras. It follows the kernel $ker(\phi)$ is a $\Dl$-ideal in $\tL_1$. Let $\tilde{I}:=I/I^2 \subseteq \P$ be the module of differentials
and let $\tL$ be any $\Dl$-Lie algebra. It follows $\tilde{I}\tL \subseteq \tL$ is an ideal. Let $(\tL, \ta,\tp,[,],D)$ be a $\Dl$-Lie algebra and let $J:=\{aD:\text{ with $a\in A$} \}\subseteq \tL$. It follows
$J$ is an ideal. Let $L:=\tL/J$. There is an exact sequence of $\P$-modules and $k$-Lie algebras
\begin{align}
&\label{ideal} 0 \rightarrow J \rightarrow \tL \rightarrow L \rightarrow 0
\end{align}
There is a canonical map of left $A$-modules and $k$-Lie algebras $\alpha_L : L \rightarrow \Der_k(A)$ making $(L, \alpha_L)$ into an $A/k$-Lie-Rinehart algebra. The sequence
\ref{ideal} is an exact sequence of $A/k$-Lie-Rinehart algebras. The ideal $J$ is a free left $A$-module on the element $D$.
\end{lemma}
\begin{proof} One checks that $ker(\phi)$ is a $\P$-submodule of $\tL_1$  and that $[\tL_1,ker(\phi)]\subseteq ker(\phi)$. Since $D\notin ker(\phi)$ it follows $ker(\phi)$ is a $\Dl$-ideal.
The same holds for $\tilde{I}$: $\tilde{I}$ is a $\P$-module and $[\tL, \tilde{I}]\subseteq \tilde{I}$. For any element $u\in \tL$ and $c\in A$ it follows 
$uc= cu+\tp(u)(c)D$. It follows 
\[ [u, aD]=a[u,D]+\tp(u)(a)D=\tp(u)(a)D \]
since $D$ is central. It follows $[\tL, J]\subseteq J$. Hence $[J,J]\subseteq J$ and $J$ is a $k$-Lie algebra. We get
\[ dc.u=uc-cu=\tp(u)(c)D \]
hence
\[ dc.dc'.u=dc(\tp(u)(c')D)=0\]
since $Dc=cD$. Hence $J$ is a left $\P$-module. It follows the Sequence \ref{ideal} is an exact sequence of $\P$-modules and $k$-Lie algebras. If $\overline{u}\in \tL/J$ is an element
with $u\in \tL$ it follows $\overline{u}c-c\overline{u}=\tp(u)(c)D:=0$. Hence $\overline{u}c=c\overline{u}$ in $\tL/J$. Hence $\tL/J$ is trivially a left $\P$-module.
Since $\tp(aD)=a\tp(D)=0$ it follows we get a canonical map
\[ \pi_L: \tL/J \rightarrow \Der_k(A) \]
and one checks $(\tL/J,\pi_L)$ is an $A/k$-Lie-Rinehart algebra. The rest is clear and the Lemma follows.
\end{proof}

\begin{definition} Let $(\tL, \ta, \tp,[,],D)$ be a $\Dl$-Lie algebra. The quotient $A/k$-Lie-Rinehart algebra $(L,\pi_L)$ with $L:=\tL/J$ from Lemma
\ref{lemma1} is the \emph{canonical quotient} of $\tL$.
\end{definition}

\begin{example} When the canonical quotient is projective. \end{example}

Assume $L:=\tL/J$ where $J$ is the ideal defined in Lemma \ref{lemma1} and let $\pi_L: L\rightarrow \Der_k(A)$ be the anchor map. 
Assume $L$ is projective as left $A$-module and let $s$ be a left $A$-linear section of the canonical projection map
$p:\tL \rightarrow L$. Define for $u,v\in L$ the following map
\[ \psi_s: \wedge^2 L \rightarrow J \]
by
\begin{align}
&\label{psi} \psi(\ou,\ov):=[s(\ou),s(\ov)]-s([\ou,\ov]) .
\end{align}
Define the map $\nabla_s: L\rightarrow \End_k(J)$ by
\[ \nabla_s(\ou)(x):=[s(\ou),x].\]

\begin{lemma} \label{Asplitting}
It follows $(J,\nabla_s)$ is a flat $L$-connection and $\psi_{s}\in \Z^2(L,(J,\nabla_s))$ where $\C^p(L,(J, \nabla_s))$ is the Lie-Rinehart complex of the flat connection $(J,\nabla_s)$.
If $s'$ is another left $A$-linear splitting of $p$ it follows there is an equality of connections $\nabla_s=\nabla_{s'}$. If $\psi_{s'}$ is the 2-cocycle associated to $s'$ it follows there is an element
$\rho\in \C^1(L,(J,\nabla_s))$ with $\psi_{s'} =\psi_s +d^1_{\nabla_s}(\rho)$. Hence there is an equality of cohomology classes 
\[ \overline{\psi_s}=\overline{\psi_{s'} }\in \H^2(L,(J,\nabla_s)).\]
\end{lemma}
\begin{proof} One checks that $\nabla_s$ is a flat $L$-connection on $J$ and that $\psi\in \Z^2(L,(J,\nabla_s))$. Assume $s'=s+\rho$ where $\rho \in \Hom_A(L,J)$. We get
for any element $u\in L$ and $x\in J$ the following:
\[ \nabla_{s'}(u)(x):=[s(u)+\rho(u),x  ]=\nabla_s(u)(x)+[\rho(u),x]=\nabla_s(u)(x) \]
since $J$ is an abelian $k$-Lie algebra. It moreover follows
\[ \psi_{s'}(u,v)=[s(u)+\rho(u),s(v)+\rho(v)]-s([u,v])-\rho([u,v])=\]
\[ \psi_s(u,v)+d^1_{\nabla_s}(\rho)(u,v)+[\rho(u),\rho(v)] \]
\[  \psi_s(u,v)+d^1_{\nabla_s}(\rho)(u,v) \]
since $J$ is an abelian Lie algebra and hence $[\rho(u),\rho(v)]=0$ for all $u,v\in L$. The Lemma follows.
\end{proof}

Define the map
\[ \phi(-,-): L\times A \rightarrow J \]
by
\begin{align}
&\label{phi} \phi(\ou,a):=s(\ou)a-s(\ou a) .
\end{align}
Make the following definition: $J\oplus^{(\phi,\psi)}L$ is the left $A$-module $J\oplus L$ with the following right $A$-module structure: Given $z:=(u,x)\in J\oplus L$ and $c\in A$
define
\[ zc:=(x, \ou)c:=(xc+\phi(\ou,c),\ou c).\]
We get since $\ou c=c\ou$ the following:
\[ \phi(\ou,c)=s(\ou)c-s(\ou c)=s(\ou)c-s(c \ou)=s(\ou)c-cs(\ou)= \]
\[ cs(\ou)+\tp(s(\ou))(c)D-cs(\ou)=\tp(s(\ou))(c)D \in J.\]
Hence
\[ zc:=cz+\tp(s(\ou))(c)D\in J\oplus L.\]
Since $ \tp = \pi_L \circ p$ and $p \circ s =Id$ it follows $\tp(s(\ou))=\pi_L(p(s(\ou))=\pi_L(\ou)$. It follows $zc=cz+\pi_L(\ou)(c)(D,0)=cz+\pi_L(\ou)(c)\tilde{D}$ where $\tilde{D}:=(D,0)$.

Let $z:=(x,\ou), w:=(y,\ov)\in J\oplus L$ and define the following product:
\begin{align}
&\label{bracket} [z,w]:=(\nabla_s(\ou)(y)-\nabla_s(\ov)(x)+\psi(\ou,\ov),[\ou,\ov]).
\end{align}
It follows the product $[,]$ is a $k$-Lie product on $J\oplus L$.  Define the map  
\[ \rho: J\oplus^{(\nabla_s, \psi)}L \rightarrow \tL \]
by
\begin{align}
&\label{rho}\rho(x,\ou):=x+s(\ou).
\end{align}
It follows $\rho$ is an isomorphism of $\P$-modules and $k$-Lie algebras. Define the map
\[ \alpha_J:J\oplus L \rightarrow \D \]
by
\[ \alpha_J := \ta \circ \rho.\]
Define $\pi_J:J\oplus L \rightarrow \Der_k(A)$ by $\pi_J(u,x):=\pi_L(x)$. Associated to a left $A$-linear splitting $s$ of $p:\tL \rightarrow L$ we get by Lemma \ref{Asplitting}
a unique flat connection $\nabla_s:L\rightarrow \End_k(J)$, a unique cohomology class $\overline{\psi_s}\in \H^2(L,(J,\nabla_s))$ and a 5-tuple
$(J\oplus L, \alpha_J, \pi_J, [,], \tilde{D})$.

\begin{proposition} \label{directsum} Let $(\tL, \ta, \tp, [,],D)$ be a $\Dl$-Lie algebra and let $(L,\alpha)$ be the canonical quotient $A/k$-Lie-Rinehart algebra of $\tL$. Assume $L$ is projective as left $A$-module. 
The 5-tupe $(J\oplus L, \alpha_L, \pi_L,[,],\tilde{D})$ constructed above is a $\Dl$-Lie algebra and there is an isomorphism $J\oplus L \cong \tL$ of $\Dl$-Lie algebras.
\end{proposition}
\begin{proof} The proof follows from the construction and calculations above.
\end{proof}

\begin{theorem}\label{Dclassification} Let $(\tL, \ta,\tp, [,],D)$ be a $\Dl$-Lie algebra and let $(L,\alpha)$ be the canonical quotient $A/k$-Lie-Rinehart algebra of $\tL$. Assume $L$ is projective as left $A$-module. 
There is an isomorphism $\tL \cong L(\fal)$ as $\Dl$-Lie algebras where $L(\fal):=F_{\fal}(L)$ and $F_{\fal}$ is the functor from Example \ref{nontrivial}. Hence $\tL$ is uniquely determined by the canonical quotient $(L,\alpha)$ and the 2-cocycle $f\in \H^2(\Der_k(A),A)$.
\end{theorem}
\begin{proof} Let $s$ be a left $A$-linear splitting of the canonical projection map $p: \tL \rightarrow L$.
From Proposition \ref{directsum} it follows the 5-tuple $(J\oplus L, \alpha_L, \pi_L,[,],\tilde{D})$ is a $\Dl$-Lie algebra and 
there is by \ref{rho} construction an isomorphism of $\Dl$-Lie algebras
\[ \rho: J\oplus L \rightarrow \tL \]
defined by
\[ \rho(aD, \ou):= aD+s(\ou) \in \tL.\]

Consider the map  $\ta: \tL \rightarrow \D$. It looks as follows: $\ta(u)=\alpha_1(u)I+\tp(u) \in A\oplus \Der_k(A)$ with $\alpha_1 \in \Hom_A(\tL, A)$

Let $z:=(aD,\ou)$ and $z':=(bD,\ov)$. It follows 
\[  [\rho(z),\rho(z')]= \]
\[ (\alpha(\ou)(b)-\alpha(\ov)(a) +\fal(\ou,\ov)+\alpha(\ou)(\alpha_1(s(\ov))-\alpha(\ov)(\alpha_1(s(\ou)))I + [\alpha(\ou), \alpha(\ov)] .\]
The Lie product on $J\oplus L$ is defined as follows:
\[ [(aD,\ou),bD,\ov)]:= (\alpha(\ou)(b)-\alpha(\ov)(a)+g(\ou,\ov))D+[\ou,\ov] \in J\oplus L\]
with 
\[ \psi(\ou,\ov):=[s(\ou),s(\ov)]-s([\ou,\ov])=g(\ou,\ov)D \in J\]
from equation \ref{psi}.

Here $g(\ou,\ov)\in \Z^2(L, A)$ is a 2-cocycle.
We get
\[ \rho([z,z'])=\]
\[ (\alpha(\ou)(b)-\alpha(\ov)(a)+g(\ou,\ov) + \alpha_1(s([\ou,\ov])) ) I +\alpha([\ou,\ov]) \in \D.\]
It follows
\[ \rho([z,z'])=[\rho(z),\rho(z')] \]
if and only if
\[  \fal(\ou,\ov)+\alpha(\ou)(\alpha_1\circ s(\ov))-\alpha(\ov)(\alpha_1 \circ s(\ou)) =g(\ou,\ov)+\alpha_1\circ s([\ou,\ov]) \]
hence $\rho$ is a map of $k$-Lie algebras if and only if

\[ g(\ou,\ov)=\fal(\ou,\ov)+d^1_{\alpha}(\alpha_1 \circ s)(\ou \wedge\ov) \]
and $\alpha_1 \circ s \in \Hom_A(L,A)$.
Hence $\rho$ is a map of $k$-Lie algebras if and only if there is an  isomorphism of $\Dl$-Lie algebras
\[ L(\fal)\cong J\oplus L .\]
It follows there is an isomorphism $\tL \cong L(\fal)$ of $\Dl$-Lie algebras and the Theorem follows.
\end{proof}

\begin{corollary} \label{mapclasscorr} Let $(\tL_i, \ta_i, \tp_i,[,],D_i)$ be $\Dl$-Lie algebras for $i=1,2$ with projective canonical quotients $(L_i,\alpha_i)$ for $i=1,2$.
There is an equality between the set of maps of $\Dl$-Lie algebras $\phi: \tL_1 \rightarrow \tL_2$ and the set of maps of $A/k$-Lie-Rinehart algebras
$\phi^*:L_1 \rightarrow L_2$.
\end{corollary}
\begin{proof} By Theorem \ref{Dclassification} there are isomorphisms $\tL_i\cong L(f^{\alpha_i})$ for $i=1,2$. It follows again from Theorem \ref{Dclassification}
that the set of maps of $\Dl$-Lie algebras between $\tL_1$ and $\tL_2$ equal the set of maps of $A/k$-Lie-Rinehart algebras between $L_1$ and $L_2$ and the Corollary 
is proved.
\end{proof}

\begin{example}\label{DLieAK} The $\Dl$-Lie algebra associated to an $A/k$-Lie-Rinehart algebra.\end{example}

Let $(L,\alpha)$ be an $A/k$-Lie-Rinehart algebra and let $\fal \in \Z^2(L,A)$ be the 2-cocycle associated  to a 2-cocycle $f\in \Z^2(\Der_k(A),A)$. In \cite{maa1}, Theorem 2.8 we constructed a functor

\[  F : \LR \rightarrow \ATiA\]
by
\[ F (L,\alpha):=(L(\fal), \alpha_f, \pi_f,[,],z) \]
where $L(\fal):=Az \oplus L$ is the abelian extension of $L$ by the free rank one $A$-module on the symbol $z$.
Define for any $\Dl$-Lie algebra $(\tL,\ta, \tp,[,],D)$ the following: $G(\tL, \ta, \tp,[,],D):=(L, \pi_L)$ where $(L,\pi_L)$ is the canonical quotient of $\tL$. Since any map of $\Dl$-Lie algebras
$\phi: \tL \rightarrow \tL'$ satisfies $\phi(D)=D'$ where $D' \in \tL'$ is the canonical central element, we get a canonical map of $A/k$-Lie-Rinehart algebras
\[ G(\phi):(L,\pi_L) \rightarrow (L', \pi_{L'}) \]
where $(L', \pi_{L'})$ is the canonical quotient of $\tL'$. Hence we get a functor
\[ G: \ATiA \rightarrow \LR .\]
Theorem \ref{Dclassification} says that in the case when the canonical quotient $L$ of $\tL$ is a projective $A$-module, it follows $G(F(L,\alpha)) \cong (L, \alpha)$ and $F(G(\tL)) \cong \tL$
are isomorphisms. 
%Hence when $L$ is a projective $A$-module it follows the functors $F$ and $G$ are inverses to each other up to isomorphism of functors.

\section{Classification of connections on $\Dl$-Lie algebras with projective canonical quotient}

In this section we classify connections on a $\Dl$-Lie algebra $\tL$ with projective canonical quotient $(L,\alpha)$. We prove in Theorem \ref{LconnLfalconn} there is a 2-cocycle $f\in \Z^2(\Der_k(A),A)$ and an equivalence of categories
\begin{align}
 &\label{equival}C_f: \Conn(L,\End) \cong \Conn(\tL ). 
\end{align}
We use the equivalence $C_f$ in \ref{equival} to classify arbitrary $\tL$-connections in Corollary \ref{DclassCorr}.

We also introduce the correspondence and Chow-operator of an $\tL$-connection $(E,\rho)$.

\begin{lemma} Let $E$ be a left $A$-module and let $\Diff^1(E)$ be the module of first order differential operators on $E$. It follows $\Diff^1(E)$ is a left $\P$-module and $k$-Lie algebra.
There is a map
\[ \psi: \Diff^1(A) \times A \rightarrow \End_A(E) \]
defined by
\[ \psi(\partial, a):= [\partial, aI]:= \partial \circ aI -aI \circ \partial \]
where $I$ is the identity operator. It follows $\psi(\partial, ab)=a\psi(\partial, b)+\psi(\partial, a)b$ for all $a,b\in A$ and $\partial \in \Diff^1(E)$.
\end{lemma}
\begin{proof} A differential operator $\partial \in \Diff^1(E)$ is by definition an operator $\partial \in \End_k(E)$ with $[[\partial, aI]bI]=0$ for all elements $a,b\in A$ where $I$ is the identity operator.
The module of differential operators $\Diff^1(E)$ has a left $A\otimes_k A$-module structure defined by $(a\otimes b. \partial)(e):=a\partial(be)$. It follows 
\[ da.\partial := \partial \circ aI-aI\circ \partial :=[\partial, aI] .\]
It follows $da.db.\partial:=[[\partial, bI],aI]=0$ hence for any element $w\in I^2$it follows $w\partial=0$ and it follows $\Diff^1(E)$ is a left $\P$-module.
One checks the product
\[ [\partial, \partial']:= \partial \circ \partial ' -\partial ' \circ \partial \]
defines a $k$-Lie algebra structure on $\Diff^1(E)$. The rest is trivial and the Lemma follows.
\end{proof}

\begin{definition} \label{lconnection} Let $(\tL,\ta,\tp,[,],D)$ be a $\Dl$-Lie algebra and let $E$ be a left $A$-module. An \emph{$\tL$-connection} $\rho$ is a map
\[ \rho: \tL \rightarrow \Diff^1(E) \]
of left $\P$-modules. The \emph{curvature} of a connection $(\rho,E)$ is the map
\[  R_{\rho}: \tL \times \tL \rightarrow \End_k(E) \]
defined by
\[ R_{\rho}(u,v):=[\rho(u),\rho(v)]-\rho([u,v]).\]
Given two $\tL$-connections $(E,\rho_E)$ and $(F,\rho_F)$, a \emph{map of  $\tL$-connections} is a map of left $A$-modules
\[ \phi: E \rightarrow F \]
with $ \rho_F(u)\circ \phi = \phi \circ \rho_E(u)$ for all $u \in \tL$. Let $\Conn(\tL)$ denote the category of $\tL$-connections and maps of connections. Let $\Conn(\tL, Id)$ denote the category
of $\tL$ connections $(\rho,E)$ with $\rho(D)=Id_E\in \End_A(E)$. 
\end{definition}

Note: It follows $\Conn(\tL, Id)$ is a full sub category of $\Conn(\tL)$. 

Note: A connection $\rho: \tL \rightarrow \Diff^1(E)$ in $\Conn(\tL,Id)$ is in particular an $A\otimes_k A$-linear map with $\rho(D)=Id_E$.

\begin{example} \label{correspondence} The degeneracy loci and correspondence associated to a connection.\end{example}

A connection in the sense of Definition \ref{lconnection} is a map of left and right $A$-modules
\[ \rho: \tL \rightarrow \Diff^1(E) \]
and we may associate to $\rho$ several types of correspondences.

\begin{definition} Given a connection $\rho \in \Hom_{\P}(\tL, \Diff^1(E))$. Let $I(\rho) \subseteq \P:=A\otimes_k A/I^2$ be the annihilator ideal
of the element $\rho$. Let $Z_{\P}(\rho):=V(I(\rho)) \subseteq \Spec(\P)$ be the \emph{correspondence of $\rho$}
\end{definition}

By definition $I(\rho)$ is the set of elements in $x\in \P$ with $x\rho=0$. The ideal $I(\rho)$ gives rise to an ideal $J(\rho)\subseteq A\otimes_k A$ containing the square of the diagonal $I^2$.
We get in a canonical way a correspondence $Z(\rho):=V(J(\rho)) \subseteq \Spec(A\otimes_k A):=X \times X$ where $X:=\Spec(A)$. Hence the connection $\rho$ gives in a canonical way rise to a correspondence
$Z(\rho)$ on $X$. 

The left $\P$-module $\Diff^1(E)$ is projective as left and right $A$-module  when $A$ is a regular ring of finite type over a field and $E$ a finite rank projective $A$-module. Hence when $\tL$ is projective as left and right $A$-module it follows a connection
\[ \rho: \tL \rightarrow \Diff^1(E) \]
is a map of left and right $A$-modules that are projective as left and right $A$-modules. Hence a connection is a geometric object and we may use $\rho$ to define a correspondence on $X\times X$.
For a classical connection 
\begin{align}
&\label{conn1} \nabla: E \rightarrow E\otimes_A \Omega^1_{A/k}
\end{align}

 it is not immediate how to do this, since $\nabla$ is a map of $k$-vector spaces and not $A$-modules. We may view
an ordinary connection as an $A$-linear map
\begin{align}
&\label{conn2} \nabla: L\rightarrow \End_k(A) 
\end{align}
where $L$ is an $A/k$-Lie-Rinehart algebra satisfying $\nabla(x)(ae)=a\nabla(x)(e)+ x(a)e$ and it is not immediate how to define a correspondence from $\nabla$.

We may associate to $\rho$ the degeneracy locies $D_{\sigma}^l(\rho)$ and $D_{\sigma}^r(\rho)$ where $l$ and $r$ refer to the degeneracy loci of $\rho$ as a left and right $A$-linear map.
Here $D_{\sigma}^l(\rho)$ and $D_{\sigma}^r(\rho)$ are defined using local trivializations of the map $\rho$. Locally the map $\rho$ is  a matrix $M$ with coefficients in a commutative ring $B$
and we may use minors of $M$ of a given size to define an ideal in $B$ associated to $\sigma$ as done in \cite{fulton}. It is possible to do this in a way that is intrinsic and does not depend 
of the choice of local trivialization  of the map $\rho$. We get
for any two $\sigma_1,\sigma_2$ a correspondence $Z(\rho, \sigma_1, \sigma_2):= D_{\sigma}^l(\rho) \times D_{\sigma}^r(\rho) \subseteq X\times X$.
Hence we may associate different types of correspondences to the connection $\rho$. If $X$ is smooth over $k$ we get for each connection $\rho$ and each correspondence $Z(\rho, \sigma_1,\sigma_2)$ 
an endomorphism
\[  I(\rho, \sigma_1, \sigma_2): \CH^*(X) \rightarrow \CH^*(X) \]
defined by
\[ I(\rho, \sigma_1,\sigma_2)(\alpha):=p_*(Z(\rho, \sigma_1,\sigma_2) \cap q^*(\alpha)) \]
where $p,q: X\times X \rightarrow X$ are the projection maps, $\CH^*(X)$ is the Chow group of  $X$ and $Z(\rho, \sigma_1,\sigma_2)\cap \alpha) $ is the intersection product. We get a similar
construction for any reasonable cohomology theory $\H(-)$ equipped with a cycle map. One would like to relate the correspondence $Z(\rho,\sigma_1,\sigma_2)$ and operation $I(\rho,\sigma_1,\sigma_2)$
to the Chern classes of $(E,\rho)$. Assume $\gamma: \CH^*(X\times_k X) \rightarrow \H^*(X\times_k X)$ is a cycle map and let $\alpha\in \H^*(X)$ be a cohomology class.
we get an operator
\[ I_{\H}(\rho, \sigma_1, \sigma_2): \H^*(X) \rightarrow \H^*(X) \]
defined by
\[  I_{\H}(\rho, \sigma_1, \sigma_2)(\alpha):=p_*( \gamma(Z(\rho, \sigma_1, \sigma_2)) \cap q^*(\alpha)) .\]

\begin{example} Algebraic cycles and the Gauss-Manin connection.\end{example}

If $\H^*(-)$ is a Weil cohomology theory, there are \emph{Lefschetz operators}
\[ \Lambda:=(L^{n-i+2})^{-1} \circ L \circ(L^{n-i}): \H^i(X) \rightarrow \H^{i-2}(X) \]
defined for any $i=0,\ldots, dim(X)$ and the operator $\Lambda$ is conjectured to be induced by an algebraic cycle $Z\subseteq X\times_k X$. The cycle $[Z]$ of the sub-scheme $Z$ induce an operator
\[ I_{\H}(Z): \H^*(X) \rightarrow \H^*(X) \]
defined by
\[ I_{\H}(Z)(\alpha):=p_*(\gamma([Z])\cap q^*(\alpha)).\]

It has been conjectured that when $X\subseteq \mathbb{P}^n_k$ is a smooth projective variety over an algebraically closed field, the Lefschetz operator $\Lambda$ is induced
by an operator on the form $I_{\H}(Z)$ for some closed sub-scheme $Z\subseteq X \times_k X$. One wants to construct non-trivial cycle classes  $\beta\in \CH^*(X\times_k X)$ and calculate the operator $I_{\H}(\beta)$.
The Chow group $\CH^*(X\times_k X)$ is hard to calculate and there are no general formulas for it. One also wants a construction of all Weil-cohomology theories $\H(-)$. If one could realize a Weil cohomology theory 
 $\H^*(-)$ as the cohomology $\H^*(\tL, -)$ of a $\Dl$-Lie algebra $\tL$ as defined in \cite{maa1}, Definition 3.23, one could approach conjectures on algebraic cycles for smooth projective families of varieties. 
One has to develop the formalism of the Gauss-Manin connection for the cohomology theory $\H^*(\tL,-)$ in the language of $\Dl$-Lie algebras. In previous papers (see \cite{maa0}) I have delevoped a formalism aimed at making explicit calculations of such connections. If one could realize the action of $\Lambda$ on a Weil cohomology $\H^*(X)$ of the total space $X$ of a smooth projective family $\pi: X\rightarrow S$, as the 
action of $\nabla(x)$ where $x$ is a vector field on $S$ and $\nabla$ the Gauss-Manin connection, this could be a first step in the direction of determining if $\Lambda$ is induced by an algebraic cycle.
The vector field $x$ and the Gauss-Manin connection $\nabla$ are algebraic objects, hence we would get an algebraic construction of $\Lambda$.

\begin{definition} \label{correspond} Let $Z(\rho,\sigma_1,\sigma_2)$ be the \emph{correspondence of $\rho$ of type $(\sigma_1,\sigma_2)$}. Let $I(\rho,\sigma_1,\sigma_2)$ be the 
\emph{Chow-operator of $\rho$ of type  $(\sigma_1,\sigma_2)$}.
\end{definition}

When using the notion of a $\Dl$-Lie algebra $\tL$, the derivation property
of the connection $\rho:\tL \rightarrow \Diff^1(E)$ is encoded in the right $A$-linearity of the map $\rho$. Hence the correspondence $Z(\rho, \sigma_1,\sigma_2)$ encodes properties of the left $A$-linearity of the map
$\rho$ and the derivation property (the right $A$-linearity) of the map $\rho$. Hence the correspondence $Z(\rho,\sigma_1,\sigma_2)$ and the endomorphism $I(\rho,\sigma_1,\sigma_2)$ 
dependes on the connection $\rho$ is a non-trivial way. It is not clear how to make a similar definition with connections on the form \ref{conn1} and \ref{conn2} depending in a non-trivial way on the derivation property
of the connection. Hence in the case of an ordinary connection or a connection $(E,\nabla)$ on an $A/k$-Lie-Rinehart algebra $L$ it is essential we work with the associated $\Dl$-Lie algebra $F_f(L):=L(\fal)$ and 
$L(\fal)$-connection $C_f(E,\nabla)$ for some 2-cocycle $f\in \Z^2(\Der_k(A),A)$ if we want to define the correspondence and Chow-operator of $(E,\nabla)$.

\begin{example} Restricted Lie-Rinehart algebras, logarithmic derivaties and Picard groups.\end{example}

Let $B$ be a commutative ring over a field $k$  characteristic $p>0$ and let $(L, \alpha)$ be a restricted $B/k$-Lie-Rinehart algebra. Let $A:=ker(L)$ be the kernel of $L$
in $B$. The ring $B$ is a \emph{purely inseparable Galois extension of $A$} if $B$ is a finitely generated and projective $A$-module and $B[L]=\Hom_A(B,B)$ in the sense of \cite{yuan}.
Yuan proves in \cite{yuan}  the existence of an exact sequence
\begin{align}
&\label{restr}  0 \rightarrow \Log(B/A) \rightarrow \Pic(A) \rightarrow \Pic(B) \rightarrow \H^2_{res}(L,B) \rightarrow \Br(B/A) \rightarrow 0 
\end{align}
where $\Log(B/A)$ is the \emph{logarithmic derivative group} of $B/A$,  $\H^2_{res}(L,B)$ is the restricted Lie-Rinehart cohomology group of $L$ with values in $B$,
$\Br(B/A)$ the Brauer group of $B/A$ and $\Pic(B), \Pic(A)$ the Picard groups of $B$ and $A$. Hence in characteristic $p>0$, the restricted version of classical $B/k$-Lie-Rinehart cohomology calculates
Picard groups and Brauer groups. One wants to generalize the exact sequence \ref{restr}  to the case of \emph{restricted $\Dl$-Lie algebras}.

\begin{example} $\Dl$-Lie algebras, groupoid schemes and algebraic stacks. \end{example}

Recall the following from \cite{ekedahl}, page 140: Let $T:=(X,Y,s_1,s_2,t,p,i)$ be a groupoid scheme with $X$ the scheme of arrows and $Y$ the scheme of objects. 
We say $T$ is a \emph{schematic equivalence relation} if the map $(s_1,s_2):X \times X \rightarrow Y \times Y$ is a \emph{monomorphism}. There is the following result:

\begin{proposition} \label{invol} Assume $Y$ is smooth over a field $k$ of characteristic zero. There is a one-to-one correspondence between locally trivial sheaves of $\O_Y/k$-Lie-Rinehart algebras
$(\mathcal{L},\alpha)$ and formal infinitesimal groupoid schemes with $Y$ as a scheme of objects. Ynder this correspondence schematic equivalence relations corresponds 
to sub-sheaves of $\O_Y/k$-Lie-Rinehart algebras of the tangent sheaf $T_{Y/k}$.
\end{proposition}
\begin{proof} See \cite{illusie}, VIII.1.1.5.
\end{proof}

Hence there is a close connection between $A/k$-Lie-Rinehart algebras, moduli spaces and algebraic stacks. If $\mathcal{F} \subseteq T_{Y/k}$ is a locally trivial finite rank sub $\O_Y$-module and 
sheaf of $k$-Lie algebras,  and the corresponding equivalence relation $R \subseteq Y\times_k Y$ is an \emph{etale equivalence relation}, we may form the \emph{stack quotient} $[Y/\mathcal{F}]$. The quotient
$[Y/\mathcal{F}]$ is a Deligne-Mumford stack. The sheaf of universal enveloping algebras $\mathcal{U}(\O_Y, \mathcal{F})$ is in a natural way a sheaf of rings on $[Y/\mathcal{F}]$ - the structure sheaf $\mathcal{O}_{[Y/\mathcal{F}]}:=\mathcal{U}(\O_Y, \mathcal{F})$  of the stack $[Y/\mathcal{F}]$. A flat $\mathcal{F}$-connection $(\mathcal{E},\nabla)$ which is a quasi coherent sheaf of $\O_Y$-modules, is in a natural way a quasi coherent sheaf of $\O_{[Y/\mathcal{F}]}$-modules.
Any Deligne-Mumford stack $[Y/R]$ arise by Proposition \ref{invol} from an involutive sub-bundle $\mathcal{F} \subseteq T_{Y/k}$. Hence if we are given the scheme of objects $Y$ of a DM-stack we may construct the equivalence relation $R$ on $Y$ using a bundle on the form of $\mathcal{F}$. The structure sheaf $\O_{[Y/R]}$ is a sheaf of filtered almost commutative rings, hence the "ringed space"
$([Y/R], \O_{[Y/R]})$ may be viewed as a non-commutative ringed space. Hence connections arise naturally when studying algebraic stacks, quasi coherent sheaves on algebraic stacks and 
non-commutative ringed spaces. 

Note: Associated to any $\Dl$-Lie algebra $\tL$ which is of finite rank and projective as left  $A$-module, we get a formal groupoid scheme over $Y:=\Spec(A)$ using the underlying $A/k$-Lie-Rinehart algebra of $\tL$
and Proposition \ref{invol}.

\begin{example} $\L$-functions and cristalline cohomology.\end{example}

Let $S:=\Spec(\O_K)$ where $K$ is an algebraic number field and let $f:X\rightarrow S$ be a regular scheme of finite type over $S$of dimension $d$. There is an equality of $\L$-functions
\begin{align}
&\label{lfunct} \L(X,s)= \prod_{p \neq (0)}\L(X(p),s) 
\end{align}
where $\L(X,s)$ is the global $\L$-function of the scheme $X$ and $\L(X(p),s)$ is the $\L$-function of the fiber $X(p):=f^{-1}(p)$ for a closed point $p \in S$. The residue field $\kappa(p)$ is a finite field
of characteristic $q >0$ and the fiber $X(p)$ is a scheme of finite type over $\kappa(p)$. Hence there is an equality
\begin{align}
&\label{zeta} \L(X(p),s)=Z(X(p), (1/q)^s) 
\end{align}
where $Z(X(p), t)$ is the Weil zeta function of $X(p)$. By the work of Kedlaya (see \cite{kedlaya}) it follows the Weil conjectures can be proved using a p-adic cohomology theory $\H^*(-)$. In fact Kedlaya has proved the Weil conjectures for an arbitrary scheme $X$ of finite type over a finite field with no condition on smoothness or projectivity on $X$. If the Ext-group $\Ext^*(V,W)$ of two connections $V,W$ calculate cristalline cohomology (see \cite{illusie}), it may be we can use $\Ext^*(V,W)$ to prove the Weil conjectures for any scheme  $X$ of finite type over a finite field. If the Weil zeta function $Z(X(p),t)$ can be calculated using the Ext group $\Ext^*(V,W)$, it may be we can use such a description in the study of the global $\L$-function $\L(X,s)$ via the product formula in \ref{lfunct}. The Ext-group is defined in complete generality and may be defined for connections on the family $X$. One has to calculate explicit examples to check if this idea leads to interesting constructions and results. If $X$ is a regular scheme the following is conjectured in \cite{soule}:
\begin{align}
&\label{sconj}  \chi(X,j)=ord_{s=j}(\L(X,s)).
\end{align}
Here
\[ \chi(X,j):= \sum_{m\geq 0} (-1)^{m+1}dim_{\mathbb{Q}}(\K_m(X)_{\mathbb{Q}}^{(d-j)} \]
is the $\K$-theoretic Euler characteristic of $X$ as defined in \cite{soule}. The conjecture in \ref{sconj} is referred to as a conjecture due to Lichtenbaum, Deligne, Bloch, Beilinson and others
in Wiles official problem description \cite{wiles} for the BSD-conjecture. If it can be proved the Ext-group $\Ext^*(V,W)$ calculate cristalline cohomology of the fibers $X(p)$ for all primes $p\neq (0)$, it might be the group $\Ext^*(V,W)$ can be used in the study of Conjecture \ref{sconj}.

%If $H \subseteq G$ is a closed subgroup of a linear algebraic group $G$ 
%of finite type over a field $k$ of characteristic zero, we may construct the quotient $\pi: G \rightarrow G/H$ and $G/H$ is a smooth quasi projective scheme of finite type over $k$.
%Any $H$-module $\rho: H \rightarrow \operatorname{GL}_k(V)$ where $V$ is a finite dimensional $k$-vector space gives rise to a finite rank locally trivial $\mathcal{O}_{G/H}$-module
%$E(\rho)$. 

\begin{lemma}\label{lemmaconn}  The following holds for an $\tL$-connection $\rho: \tL \rightarrow \Diff^1(E)$: 
\[ \rho(u)(ae)=a\rho(u)(e)+\tp(u)(a)\rho(D)(e).\]
for all elements $a\in A, u\in \tL$ and $e\in E$. It follows $\rho(D)\in \End_A(E)$. The curvature $R_{\rho}$ defines a map
\[ R_{\rho}:\tL \times \tL \rightarrow \Diff^1(E).\]
Assume $P:\tL \rightarrow \End_A(E)$ is an $\P$-linear map. It follows $\rho':=\rho+P$ is an $\tL$-connection. If $P(D)=0$ it follows $\rho'(D)=Id_E$.
\end{lemma}
\begin{proof} Since $\rho$ is $\P$-linear it follows $\rho$ is $A\otimes_k A$-linear. It follows $\rho(au)=a\rho(u)$ and $\rho(ua)=\rho(u)a$. We get
\[ \rho(u)(ae)=\rho(ua)(e)=\rho(au+\tp(u)(a)D)(e)=a\rho(u)(e)+\tp(u)(a)\rho(D).\]
Since $Da=aD$ it follows $\rho(D)\in \End_A(E)$. The second statement holds since the element $[\rho(u),\rho(v)]:=\rho(u)\rho(v)-\rho(v)\rho(u)\in \Diff^1(E)$. 
Since $\rho$ is $\P$-linear it follows $\rho':=\rho + P$ is $\P$-linear. If $P(D)=0$ it follows $\rho'(D)=\rho(D)+P(D)=Id_E$.
The Lemma follows.
\end{proof}

Hence the notion of an $\tL$-connection introduced in Definition \ref{lconnection} agrees with the notion introduced in the paper \cite{maa1}.

\begin{example} The $L(\fal)$-connection associated to an $(L,\psi)$-connection.\end{example}

\begin{definition}
Let $(L,\alpha)$ an $A/k$-Lie-Rinehart algebra and let $\nabla:L\rightarrow \End_k(E)$ be an $(L,\psi)$-connection where $\psi\in \End_A(E)$. This means 
\[ \nabla(x)(ae)=a\nabla(x)(e)+\alpha(x)(a)\psi(e) \]
for all $a\in A, e\in E$ and $x\in L$. Let $\Conn(L,\End)$ denote the category of $(L,\psi)$-connections and morphisms. The endomorphism
$\psi$ may vary. Let $\Conn(L)$ denote the category of ordinary $L$-connections and morphisms of $L$-connections.
%Let $\Mod(L(\fal))$ denote the category of $L(\fal)$-connections and morphisms of $L(\fal)$-connections. Let $\Mod(L(\fal),Id)$ denote the category of $L(\fal)$-connections
%\[ \rho: L(\fal) \rightarrow \Diff^1(E) \]
%with $\rho(z)=Id$.
\end{definition}

Note: If $\nabla:L\rightarrow \End_k(E)$ is an ordinary connection and $\psi\in \End_A(E)$ it follows $\nabla \circ \psi$ is an $(L,\psi)$-connection.

Recall the following construction:
Let $f\in \Z^2(\Der_k(A),A)$ be a 2-cocycle and let $F_f(L):=(L(\fal), \alpha_f,\pi_f,[,],z)$ be the $\Dl$-Lie algebra associated to $L$ and $f$. Define the following map:
\[ \rho:=\rho_{\nabla}:L(\fal)\rightarrow \End_k(E) \]
by
\[ \rho(az+x):=a\psi+\nabla(x) \in \End_k(E) .\]
Let $C_f(E,\nabla):=(E,\rho_{\nabla})$.

Let $(F,\nabla')$ be an  $(L,\psi')$-connection with $\psi'\in \End_A(E)$ and let $\phi:(E,\nabla) \rightarrow (F,\nabla')$ be a map of connections.
Define the following map
\[  C_f(\phi):=\phi:E\rightarrow F .\]
It follows the map $C_f(\phi):(E,\rho_{\nabla})\rightarrow (F,\rho_{\nabla'})$ is a map of $L(\fal)$-connections.
Let $i:L \rightarrow L(\fal)$ be the canonical left $A$-linear map and let $(E,\rho)$ be an $L(\fal)$-connection:
\[ \rho: L(\fal) \rightarrow \Diff^1(E) .\]
 Let $\nabla_{\rho}:=\rho \circ i$. It follows
\[ \nabla:L \rightarrow \End_k(E) \]
is an $(L,\psi)$-connection. Let $R_f(E,\rho):=(E, \nabla_{\rho})$.

\begin{lemma} \label{mlemma} Let $(L,\alpha)$ be an $A/k$-Lie-Rinehart algebra and let $(E,\nabla)$ be an $L$-connection. Let $(F,\rho)$ be an $L(\fal)$-connection
The construction above define functors 
\[   C_f: \Conn(L,\End) \rightarrow \Conn(L(\fal)) \]
by
\[ C_f(E,\nabla):=(E,\rho_{\nabla}) \]
and
\[ R_f: \Conn(L(\fal)) \rightarrow \Conn(L,\End) \]
by
\[ R_f(F,\rho):=(F, \nabla_{\rho}).\]
It follows $C_f \circ R_L = Id$ and $R_f \circ C_f=Id$ hence $C_f$ and $R_f$ are equivalences of categories for any 2-cocycle $f\in \Z^2(\Der_k(A),A)$
\end{lemma}
\begin{proof} The proof follows immediately from the constructions above.
\end{proof}

\begin{theorem}\label{LconnLfalconn}  Let $(L,\alpha)$ be an $A/k$-Lie-Rinehart algebra and let $f\in \Z^2(\Der_k(A),A)$ be a 2-cocycle. Lemma \ref{mlemma} gives an equivalence of categories
\[ C_f: \Conn(L,\End) \rightarrow \Conn(L(\fal)) \]
from the category $\Conn(L,\End)$ of $(L,\psi)$-connections, to the category $\Conn(L(\fal))$ of $L(\fal)$-connections $\rho$ for any 2-cocycle $f$.
Let $u:=az+x, v:=bz+y\in L(\fal)$ and let $e\in E$. Let $(E,\nabla)$ be an $L$-connection and let $\rho_{\nabla}:=C_f(\nabla)$. The following holds:
\[ R_{\rho_{\nabla}}(u,v)(e) =R_{\nabla \circ i}(x,y)(e)-\fal(u,v)e.\]
For any $L(\fal)$-connection $(E,\rho)$ there is an $(L,\psi)$-connection $(E,\nabla)$ with $C_f(E,\nabla)=(E,\rho)$. If $\rho(z)=Id$ we may chose $(E,\nabla)\in \Conn(L,Id)$.
\end{theorem}
\begin{proof} One checks for any $(L,\psi)$-connection $(E,\nabla)$ the corresponding map $\rho_{\nabla}$ is a map
\[ \rho_{\nabla}: L(\fal) \rightarrow \Diff^1(E) \]
of left $\P$-modules. Moreover for any map of $(L,\psi)$-connections $\phi: (E,\nabla) \rightarrow (F, \nabla')$ it follows the map
\[ C_f(\phi):C_f(E,\nabla) \rightarrow C_f(F,\nabla') \]
is a map of $L(\fal)$-connections with $C_f(\phi \circ \psi)=C_f(\phi) \circ C_f(\psi)$. By definition $C_f \circ R_f = Id$ and $R_f \circ C_f = Id$ and the first claim follows.
The statement on the curvature follows from Lemma 2.20 in \cite{maa1}.
In particular given an $L(\fal)$-connection $(E,\rho)$ it follows  $(E,\rho)\cong C_f(R_f(E,\rho))$ and $R_f(E,\rho)\in \Conn(L,\End)$.
The Theorem follows.
\end{proof}

We may classify $\tL$-connections in terms of $(L,\psi)$-connections in the case when the canonical quotient $L$ of $\tL$ is a projective $A$-module.

\begin{corollary} \label{DclassCorr} Let $(\tL,\ta,\tp,[,],D)$ be a $\Dl$-Lie algebra and let $(E,\rho)$ be an $\tL$-connection. Assume the canonical quotient $L$ of $\tL$
is a projective $A$-module. It follows any $\tL$-connection $(E,\rho)$ is on the form $C_f(E,\nabla)$ where $f\in \Z^2(\Der_k(A),A)$ and $(E,\nabla)$ is an $(L,\psi)$-connection for some
$\psi \in \End_A(E)$. If $\rho(D)=I$ it follows there is an $L$-connection $(E,\nabla)$ with $C_f(E,\nabla)=(E,\rho)$.
\end{corollary}
\begin{proof} By Theorem \ref{Dclassification} there is an isomorphism $\tL\cong L(\fal)$ for $f\in \Z^2(\Der_k(A),A)$. The Corollary now follows from Theorem \ref{LconnLfalconn}.
\end{proof}

\begin{example} Non-abelian extensions of $\Dl$-Lie algebras. \end{example}

\begin{definition} \label{nonab} Let $(\tL, \ta, \tp,[,],D)$ be a $\Dl$-Lie algebra and let $(\rho,E)$ be an $\tL$-connection with $\rho(D)=Id_E$. Let
$(\End(\tL,E), \alpha_E, \pi_E, [,],\tilde{D})$ be the $\Dl$-Lie algebra constructed in Section 4 in \cite{maa1}.
\end{definition}

The following Theorem generalize Theorem 2.14 from \cite{maa2} from a connection on an $A/k$-Lie-Rinehart algebra to the case of a connection on a $\Dl$-Lie algebra:

\begin{theorem} \label{mthm} Let $(\tL, \ta,\tp,[,],D)$ be a $\Dl$-Lie algebra and let $(\rho,E)$ be an $\tL$-connection with $\rho(D)=Id_E$. Let $\End(\tL,E)$ be the non-abelian extension
of $\tL$ with $\rho$ from Definition \ref{nonab}. There is a canonical flat connection
\[ \rho^! : \End(\tL,E)\rightarrow \Diff^1(E) \]
defined by $\rho^1(\phi, u):=\phi+\rho(u)$. The exact sequence
\begin{align}
&\label{exts}  0 \rightarrow \End_A(E) \rightarrow \End(\tL, E) \rightarrow \tL \rightarrow 0
\end{align}
is split in the category of $\Dl$-Lie algebras if and only if $E$ has a flat $\tL$-connection.
\end{theorem}
\begin{proof} By Proposition 4.9 in \cite{maa1} we get an extensins of $\Dl$-Lie algebras
\[ 0 \rightarrow \End_A(E) \rightarrow \End(\tL,E) \rightarrow \tL \rightarrow 0\]
where $\End(\tL,E):= \End_A(E)\oplus \tL$ with $\tilde{D}:=(0,D)$ and 
$p_E: \End(\tL,E)\rightarrow \tL$ defined by $p_E(\phi,u):=u\in \tL$. Define $c(\phi,u):=(c\phi, cu)$ and $(\phi,u)c:=(\phi c,uc)$ and define
\[ \alpha_E: \End(\tL,E) \rightarrow \D \]
by
\[\alpha_E(\phi,u):= \ta(u)\in \D \]
and
\[ \pi_E:\End(\tL,E) \rightarrow \Der_k(A) \]
by
\[\pi_E(\phi,u):=\tp(u).\]
Define moreover for any $z:=(\phi,u),w:=(\psi,v)\in \End(\tL,E)$
\[ [z,w]:=([\phi,\psi]+[\rho(u),\psi]-[\rho(v),\phi]+R_{\rho}(u,v),[u,v]).\]
It follows from Proposition 4.9 in \cite{maa1} that $\End(\tL,E)$ is an extension of $\tL$ by the $A$-Lie algebra $\End_A(E)$. Asssume $s:\tL \rightarrow \End(\tL,E)$ is a section 
of the map $p_E$. Hence $p_E \circ s=Id_{\tL}$ and $s$ is a map of $\Dl$-Lie algebras. It follows $s$ is $\P$-linear, a map of $k$-Lie algebras and $s(D)=(0,D)$. It follows 
$s(u)=(P(u),u)$ where $P: \tL \rightarrow \End_A(E)$ is a $\P$-linear map with $P(D)=0$. It follows from Lemma \ref{lemmaconn} the map $\rho':=\rho+P$ is an $\tL$-connection
$\rho': \tL \rightarrow \Diff^1(E)$. One checks that the map $s$ is a map of $k$-Lie algebras if and only if $R_{\rho'}=0$, hence the sequence \ref{exts} is split in the category of $\Dl$-Lie algebras
if and only if $E$ has a flat $\tL$-connection. The Theorem is proved.
\end{proof}

Note: In \cite{maa2}, Theorem 2.4 a result similar to Theorem \ref{mthm} is proved for $A/k$-Lie-Rinehart algebras.  Note moreover that by Theorem \ref{mthm} we may view any $\tL$-connection
$(\rho,E)$ as a representation of the $k$-Lie algebra $\End(\tL,E)$. Since the induced connection $\rho^!$ is flat, it follows $\rho^!$ is a map of $k$-Lie algebras. 

\begin{corollary} \label{nonabclass} Let $(\tL, \ta,\tp,[,],D)$ be a $\Dl$-Lie algebra and let $(\rho,E)$ be an $\tL$-connection with $\rho(D)=Id_E$. It follows the canonical quotient of $\End(\tL,E)$ is isomorphic to the  
$A/k$-Lie-Rinehart algebra $\End(L,E)$ where $(L, \pi_L)$ is the canonical quotient of $\tL$. If $L$ is projective there is an isomorphism $\End(\tL,E)\cong \End(L(\fal), E)$, 
where $L(\fal)=F_{\fal}(L)$ and $F_{\fal}$ is the functor from Example \ref{nontrivial}.
\end{corollary}
\begin{proof} The proof follows from Theorem \ref{mthm} and Theorem \ref{Dclassification} since $\End(\tL,E):=\End_A(E)\oplus \tL$ and $\tilde{D}:=(0,D)$.
\end{proof}

\end{document}